\documentclass[11pt,a4paper]{article}

\usepackage[utf8]{inputenc}
\usepackage[T1]{fontenc}
\usepackage{lmodern}
\usepackage[american]{babel}
\usepackage{etex}

\usepackage[margin=3cm]{geometry}

\usepackage{csquotes}
\usepackage{amsmath,amsfonts,amstext,amssymb,amsbsy,amsopn,eucal,dsfont,mathtools}
\mathtoolsset{mathic}
\usepackage[amsmath,thmmarks]{ntheorem}

\usepackage{graphicx}
\usepackage[shortlabels]{enumitem}
\usepackage[vlined,ruled]{algorithm2e}
\usepackage{xspace}
\usepackage{tikz}
\usetikzlibrary{decorations.pathreplacing}
\usepackage{microtype}
\usepackage{float}

\usepackage{hyperref}
\usepackage[capitalize]{cleveref}

\newcommand{\N}{\mathbb{N}}
\newcommand{\R}{\mathbb{R}}
\newcommand{\trace}{\ensuremath{\operatorname{trace}}}

\newlist{compactenum}{enumerate}{4}
\setlist[compactenum]{label=\arabic*.,itemsep=0pt,topsep=0pt,parsep=0pt,leftmargin=*}
\newlist{compactitem}{itemize}{4}
\setlist[compactitem]{label=\textbullet,itemsep=0pt,topsep=0pt,parsep=0pt}

\crefformat{enumi}{#2#1#3}
\crefrangeformat{enumi}{#3#1#4--#5#2#6}
\crefmultiformat{enumi}{#2#1#3}{ and~#2#1#3}{, #2#1#3}{ and~#2#1#3}
\crefrangemultiformat{enumi}{#3#1#4--#5#2#6}{ and~#3#1#4--#5#2#6}{,#3#1#4--#5#2#6}{ and~#3#1#4--#5#2#6}
\crefalias{enumii}{enumi}
\crefalias{enumiii}{enumi}
\crefalias{compactenumi}{enumi}
\crefalias{compactenumii}{enumi}
\crefalias{compactenumiii}{enumi}
\crefformat{equation}{(#2#1#3)}
\crefrangeformat{equation}{(#3#1#4)--(#5#2#6)}
\crefmultiformat{equation}{(#2#1#3)}{ and~(#2#1#3)}{, (#2#1#3)}{ and~(#2#1#3)}
\crefrangemultiformat{equation}{(#3#1#4)--(#5#2#6)}{ and~(#3#1#4)--(#5#2#6)}{, (#3#1#4)--(#5#2#6)}{ and~(#3#1#4)--(#5#2#6)}

\theoremstyle{plain}
\theorembodyfont{\normalfont\itshape}
\theoremseparator{.}
\newtheorem{Theorem}{Theorem}
\newtheorem{Lemma}[Theorem]{Lemma}
\newtheorem{Observation}[Theorem]{Observation}

\newtheorem{Conjecture}[Theorem]{Conjecture}
\newtheorem{Example}[Theorem]{Example}

\theoremstyle{proof}
\theorembodyfont{\normalfont}
\theoremsymbol{\ensuremath{\square}}
\newtheorem*{Proof}{Proof}
\newcommand{\ie}{\emph{i.\,e\mbox{.}}\xspace}

\author{Christoph Helmberg\thanks{Fakult\"at f\"ur Mathematik, Technische
Universit\"at Chemnitz, D-09107 Chemnitz,    Germany. helmberg@mathematik.tu-chemnitz.de} 
\and Guilherme Porto\thanks{Instituto Federal de Educação, Ciência e Tecnologia Farroupilha, Campus São Borja. guilherme.porto@iffar.edu.br} \and Guilherme Torres\thanks{In Memory: During the preparation of this note, Guilherme Torres tragically passed away.}
\and Vilmar Trevisan\thanks{Instituto de Matem\'atica, Universidade
Federal do Rio Grande do Sul, Porto Alegre. Brazil. trevisan@mat.ufrgs.br} }%

\title{An interlacing property of\\ the signless Laplacian of threshold graphs}
\date{\today}

\begin{document}

\maketitle
\begin{abstract}

%   We show that for threshold graphs the eigenvalues of the signless
%   Laplacian interlace with the node degrees greater or equal to the
%   $i$-th smallest degree $d_{i-1}$ \cvt{$d_i$?}and less or equal to the $i-1$-st smallest
%   degree $d_{i-1}$ for the unique  $i$ satisfying $d_{i-1}\le n-i \le d_{i}$. As an application, we show that the signless Brouwer conjecture holds for threshold graphs, i.e., for   threshold graphs the sum of the $k$ largest eigenvalues is bounded by the number of edges plus $k+1$ choose 2.
% \vspace{.5cm}

We show that for threshold graphs, the eigenvalues of the signless
  Laplacian matrix interlace with the degrees of the vertices. As an application, we show that the signless Brouwer conjecture holds for threshold graphs, i.e., for   threshold graphs the sum of the $k$ largest eigenvalues is bounded by the number of edges plus $k+1$ choose 2.

  \noindent\mbox{\textbf{Keywords:} threshold graphs, signless Laplacian
  spectrum, Brouwer conjecture}


\end{abstract}

\section{Introduction and main result}
Given an undirected graph $G=(N,E)$ on node set $N=[n]:=\{1,\dots,n\}$
for some $n\in\N$ and
edge set $E\subseteq {N\choose 2}:=\{\{i,j\}\colon i,j\in N,i\neq j\}$,
the signless Laplacian is the symmetric $n\times n$ matrix
$Q(G)=\sum_{\{i,j\}\in E} (e_i+e_j)(e_i+e_j)^\top$, where $e_i$ denotes
the $i$-th column of the $n\times n$ identity matrix $I_n$. $Q(G)$ may
also be written as $Q(G)=D+A$, where $A$ is the adjacency matrix and
$D$ the diagonal degree matrix. Thus, $Q(G)$ is positive semidefinite and its eigenvalues may be ordered as $$0\leq \lambda_1 \leq \lambda_2 \leq \cdots \leq \lambda_n.$$

Spectral properties of the signless Laplacian $Q(G)$ of a graph $G$ were collected and developed in a series of papers by Cvetkovi{\'c} and Simi{\'c} \cite{cvetkovic09,cvetkovic10,cvetkovic11}. Nevertheless, it seems that the spectrum of this matrix is far less understood than that of the combinatorial Laplacian $L(G) = D - A$, for example.

In this note, we show an interlacing property of the eigenvalues of
$Q(G)$ when $G$ is a threshold graph. Essentially, our result says
that the eigenvalues of a threshold graph interlace with the degrees of the vertices.

In order to more precisely state the result, we first review some
facts about threshold graphs. This class of graphs has been discovered
independently by several authors in many distinct contexts since the
1970's. They are an important class of graphs because of their
numerous applications in diverse areas which include computer science,
social sciences and psychology. See, for example, \cite{mahadev95} for
a more detailed account. A threshold graph can be characterized in
many ways. We are going to view threshold graphs as given through an
iterative process steered by a binary string which starts with an
isolated vertex\ (for the initial digit 0 or 1), and where, at each step, either a new isolated vertex is added (digit 0), or a dominating
vertex (adjacent to all previous vertices, digit 1) is added.

Building on the notation in \cite{FritscherTrevisan2016}, this
construction is encoded in a binary sequence $b_1^{q_1}\dots
b_r^{q_r}$ of length $n=\sum_{k=1}^rq_k$ with $b_k\in\{0,1\}$ and
$q_k\in\N$ giving this digit's repetitions in the sequence for $k=1,\dots,r$.
For convenience, let $n_k=\sum_{h=1}^kq_k$ for $k=0,\dots,r$, so $n=n_r$
and $n_0=0$. Then the corresponding threshold graph $G=(N=[n],E)$
has edge set $E$ so that for $1\le i<j\le n$
\begin{equation} \label{eq:TGdef}
 ij\in E\quad \Leftrightarrow\quad \exists k\in\{1,\dots,r\} \ (n_{k-1}< j\le n_k \wedge b_k=1).
\end{equation}
All nodes of block $k$ have the same degree. Their values are
\begin{equation}
  \label{eq:pidef}
  p_k=\left\{
    \begin{array}{ll}
      \sum_{h=k}^{r}b_hq_{h} &\text{ for }b_k=0,\\
      \sum_{h=1}^{k}q_h-1+\sum_{h=k+1}^{r}b_hq_{h} &\text{ for }b_k=1,\\
    \end{array}\right.
\qquad k=1,\dots,r.
\end{equation}
We may and will assume the sequence of
$b_k$ to be alternating, \ie, $b_k+b_{k+1}=1$ for $k=1,\dots,r-1$.
Because $G$ is independent of the choice of the first digit, there is
no loss in generality in requiring $q_1\ge 2$.
It is known (see, for example,
\cite{cvetkovic2007,FritscherTrevisan2016}) that $p_k-b_k$ is an
eigenvalue of $Q(G)$ of multiplicity at least $q_k-1$, for
$k=1,\ldots,r$. As a direct consequence of the assumptions,
the degrees and corresponding eigenvalues satisfy
\begin{equation}
  \label{eq:pi}
  \begin{array}{c}
  p_{r-b_r} < p_{r-b_r-2} < \dots
  < p_{1+b_1+2} < p_{1+b_1}\le  p_{2-b_1}-1, \\[1ex]
  p_{2-b_1}-1 < p_{4-b_1}-1 < p_{6-b_1}-1 < \dots < p_{r-(1-b_r)}-1.
\end{array}
\end{equation}
The results of \cite{FritscherTrevisan2016} imply that the remaining
$r$ eigenvalues are those of
the following \emph{condensed signless Laplacian} $r\times r$ matrix
\begin{displaymath}
  C(G)=\left[
    \begin{array}{ccccc}
      p_1+b_1(q_1-1)  & b_2\sqrt{q_1q_2} & b_3\sqrt{q_1q_3} & \dots & b_r\sqrt{q_1q_r} \\
      b_2\sqrt{q_1q_2} & p_2+b_2(q_2-1)& b_3\sqrt{q_2q_3} & \dots & b_r\sqrt{q_2q_r} \\
      b_3\sqrt{q_1q_3} & b_3\sqrt{q_2q_3} & p_3+b_3(q_3-1)& \dots & b_r\sqrt{q_3q_r} \\
      \vdots & & & & \vdots\\
      b_r\sqrt{q_1q_r} & b_r\sqrt{q_2q_r} & b_r\sqrt{q_3q_r} & \dots & p_r+b_r(q_r-1)
    \end{array}
\right].
\end{displaymath}
In \cref{T:condensedeigs} we establish that these remaining signless
Laplacian eigenvalues of the threshold graph $G$ interlace with
$p_k-b_k$. More precisely, for $C=C(G)$ we prove the inequalities
  \begin{displaymath}
   \lambda_1(C)\le p_{r-b_r} \le \lambda_2(C) \le p_{r-b_r-2} \le \dots
  \le \lambda_{\frac{r-b_r-b_1+1}2}(C) \le p_{1+b_1}
\end{displaymath}
and
  \begin{displaymath}
  p_{2-b_1}-1\le \lambda_{\frac{r-b_r-b_1+1}2+1}(C)\le p_{4-b_1}-1\le\lambda_{\frac{r-b_r-b_1+1}2+2}(C)\le
  \dots \le p_{r-(1-b_r)}-1\le \lambda_r(C).
\end{displaymath}
This interlacing property gives rise to the main  result \cref{T:main}
which sheds some light on the distribution of the signless
Laplacian spectrum of graphs. The trivial upper bound $2n-1$ on the
largest signless Laplacian eigenvalue of a graph $G$ --- called the
signless spectral radius of $G$ --- is frequently an obstacle for
obtaining meaningful upper bounds involving this spectral parameter.
For threshold graphs, however, this interlacing property may be used
to obtain tighter bounds for the sum of the largest eigenvalues.
Indeed, together with a subtle sharpening for a special case, the
bounds suffice to  establish the signless Brouwer conjecture
\cite{ASHRAF20134539} for threshold graphs in \cref{T:signlessBrouwer}.

The remainder of the paper is organized as follows. The next section is devoted
to prove our main result. In \cref{sec:disc} we interpret the result
and give an estimate of the eigenvalues growth with the evolution of a
binary sequence defining a threshold graph. In \cref{sec:brouwer}, as an application of the interlacing property for threshold graphs, we
prove that the signless Brouwer conjecture holds for threshold graphs.

\section{Proof of the main result}\label{sec:main}

We will make heavy use of the following direct consequence of the Courant-Fischer Theorem.
\begin{Theorem}\label{T:CF}
Given a symmetric matrix $A\in \R^{n\times n}$ and $h\in\{1,\dots,n\}$,
for any choice of vectors $v_1,\dots,v_{h-1}\in\R^n$ and
$u_1,\dots,u_{n-h}\in\R^n$ there holds
\begin{displaymath}
  \min_{
    \begin{array}{c}
x\bot v_1,\dots,v_{h-1}\\
\|x\|=1
\end{array}} x^\top Ax \quad\le\quad \lambda_h(A)\quad \le\quad \max_{
    \begin{array}{c}
x\bot u_1,\dots,u_{n-h}\\
\|x\|=1
\end{array}}x^\top Ax.
\end{displaymath}
If $h-1$ or $n-h$ is zero, $x$ ranges over $\R^n$ in the respective
expression.
\end{Theorem}
We recall and use a number of further consequences of the
Courant-Fischer Theorem.
%\cvt{Inside the [] it doesn't accept this: \cite[Theorem 4.3.1]{hornMatrixAnalysis}.}

\begin{Theorem} [Weyl Inequalities -- \cite{hornMatrixAnalysis},Theorem 4.3.1] \label{lemma:weylIneq}
Let $M$ and $N$ be two $n \times n$ real symmetric matrices. Then
$$
\lambda_{i}(M+N) \leq \lambda_{i+j}(M) + \lambda_{n-j}(N),
$$
for $1 \leq i \leq n$ and $0 \leq j \leq n-i$, and
$$
\lambda_i(M+N) \geq \lambda_{i-j+1}(M) + \lambda_{j}(N)
$$
for $1 \leq i \leq n$ and $1 \leq j \leq i$.
\end{Theorem}
\begin{Theorem}[Cauchy Interlacing -- \cite{hornMatrixAnalysis},Theorem 4.3.28]\label{lemma:cauchy}
Let $A$ be a real, symmetric $n \times n$ matrix and let $B$ be a principal submatrix of $A$ with order $m \times m$. Then, for $k = m, \cdots, 1$,
$$
\lambda_{k+n-m}(A) \geq \lambda_{k}(B) \geq \lambda_k(A).
$$
\end{Theorem}
\begin{Theorem}[Interlacing for a rank-one perturbation --
  \cite{hornMatrixAnalysis}, Corollary
  4.3.9]\label{interlacing_rank_one}\quad\\
Let $A \in \R^{n\times n}$ be a symmetric matrix and let $z \in \R^n$ be a vector. Then the eigenvalues of $A$ and $A + zz^\top$  satisfy
\begin{eqnarray}
\begin{array}{c}
\lambda_h(A) \leq \lambda_h(A + zz^\top) \leq \lambda_{h+1}(A) , \ \ h = 1, 2, \ldots, n-1, \\[2mm]
\lambda_n(A) \leq \lambda_{n}(A + zz^\top). \\
\end{array}
\end{eqnarray}
\end{Theorem}

In order to relate the eigenvalues of the signless Laplacian matrix of a
general graph $G$ to those of its complement, note that the signless Laplacian of the complete
graph is $Q(K_n)=\sum_{\{i,j\}\in {[n]\choose 2}}
(e_i+e_j)(e_i+e_j)^\top=(n-2)I_n+\mathbf{1}_n\mathbf{1}_n^\top$,
 where $\mathbf{1}$ denotes the vector of all ones of given or appropriate
dimension. With this the signless Laplacian of the complement graph $\bar G=(N,{N\choose 2}\setminus
E)$ of $G$ computes to $Q(\bar G)=Q(K_n)-Q(G)=(n-2)I_n+\mathbf{1}_n\mathbf{1}_n^\top-Q(G)$.

\begin{Lemma}\label{L:Cauchy}
  Let $Q=Q(G)$ be the signless Laplacian of a graph $G$ on
  $n$ nodes and let $\bar Q=Q(\bar G)$ represent its complement, then
  \begin{displaymath}
    \begin{array}{r@{\;}l}
    \max\{n-2-\lambda_n(\bar Q),0\}&\le \lambda_1(Q)\le n-2-\lambda_{n-1}(\bar
    Q)\le\lambda_2(Q)\le \cdots\\
 \cdots& \le \lambda_{n-1}(Q)\le n-2-\lambda_1(\bar Q)\le \min\{n-2,\lambda_n(Q)\}.
\end{array}
  \end{displaymath}
\end{Lemma}
\begin{Proof}
  The signless Laplacian is positive semidefinite,
  so $\lambda_1(Q)\ge 0$ and $\lambda_1(\bar Q)\ge 0$. By
  $Q=\mathbf{1}\mathbf{1}^\top+(n-2)I_n-\bar Q$ the
  statement follows from \cref{interlacing_rank_one}.
\end{Proof}

If now $G$ is a threshold graph, we notice that Equ.\ \eqref{eq:TGdef} establishes that the complement graph $\bar G$ of $G$ is defined by the binary sequence $(1-b_1)^{q_1}\dots(1-b_r)^{q_r}$. This sequence is again alternating with $q_1\ge 2$. The degrees $\bar p_i$ of the blocks of the
complement graph satisfy
  \begin{equation}
    \label{eq:compldeg}
    \bar p_k=n-1- p_k \quad\text{ for } k=1,\dots,r.
  \end{equation}

Up to the specification of the eigenvectors the following result has already been observed in \cite{FritscherTrevisan2016}.
\begin{Lemma}\label{L:directvecs}
 Let $b_1^{q_1}\dots b_r^{q_r}$ be a binary sequence specifying
  a threshold graph $G$ on $n=\sum_{k=1}^rq_k$ nodes with block degrees
  $p_k$  and put $n_0=0$, $n_k=\sum_{h=1}^kq_h$ for $k=1,\dots,r$. For
  $k\in\{1,\dots,r\}$ its signless Laplacian
  $Q(G)$ has an eigenvalue $p_k-b_k$ of multiplicity at least $q_k-1$ with
  an associated orthogonal basis of eigenvectors
  \begin{displaymath}
    v^{(k)}_1=\left[
      \begin{array}{c}
        \mathbf{0}_{n_{k-1}}\\1\\-1\\\mathbf{0}_{n-n_{k-1}-2}
      \end{array}
\right],v^{(k)}_2=\left[
      \begin{array}{c}
        \mathbf{0}_{n_{k-1}}\\\frac12\mathbf{1}_2\\-1\\\mathbf{0}_{n-n_{k-1}-3}
      \end{array}
\right],
\dots,v^{(k)}_{q_k-1}=\left[
      \begin{array}{c}
        \mathbf{0}_{n_{k-1}}\\\frac1{q_k-1}\mathbf{1}_{q_k-1}\\-1\\\mathbf{0}_{n-n_{k-1}-q_k}
      \end{array}
\right].
  \end{displaymath}
\end{Lemma}
\begin{Proof}
  Direct computation shows $(v^{(k)}_i)^\top v^{(k)}_j=0$ for $1\le i<
  j<q_k$. Let $v$ be any of these vectors.
  By \eqref{eq:TGdef} columns $n_{k-1}+1$ to $n_k$ of row $i\le n_{k-1}$ or row
  $i>n_{k}$ of $Q$ have the same value $b_k$, so $(Qv)_i=0$ for these
  $i$. It remains to consider the principal submatrix $Q'$ on
  indices $i=n_{k-1}+1,\dots,n_k$. The case $b_k=0$ yields
    $Q'=p_kI_k$, so  $Qv=p_kv$; the case $b_k=1$ yields
    $Q'=\mathbf{1}_{q_k}\mathbf{1}_{q_k}^\top+(p_k-1)I_{q_k}$, so
    $Qv=(p_k-1)v$.
\end{Proof}
\cref{L:directvecs} describes $n-r$ eigenvalues of $Q$ via its
eigenvectors. It is proved in \cite{FritscherTrevisan2016} that the
remaining $r$ eigenvalues are those of $C=C(G)$.
By \eqref{eq:compldeg}
the condensed signless Laplacian $\bar C=C(\bar G)$ of the complement graph may
be computed via
\begin{displaymath}
\bar C=(n-2)I_r+\hat q\hat q^\top-C\quad\text{ with } \hat
q=(\sqrt{q_1},\dots,\sqrt{q_r})^\top.
\end{displaymath}
So the eigenvalues of $C$ and $\bar C$ satisfy the same interlacing
property as those of $Q$ and
$\bar Q$.
\begin{Lemma}\label{L:condensedCauchy}
  Let $C$ be the condensed signless Laplacian of a threshold graph $G$
  specified by a binary sequence $b_1^{q_1}\dots b_r^{q_r}$ on
  $n=\sum_{k=1}^rq_k$ nodes and let $\bar C$ represent its corresponding complement, then
  \begin{displaymath}
    \begin{array}{r@{\;}l}
    \max\{n-2-\lambda_r(\bar C),0\}&\le \lambda_1(C)\le n-2-\lambda_{r-1}(\bar
    C)\le\lambda_2(C)\le \cdots\\
 \cdots& \le \lambda_{r-1}(C)\le n-2-\lambda_1(\bar C)\le \min\{n-2,\lambda_r(C)\}.
\end{array}
  \end{displaymath}
\end{Lemma}
\begin{Proof}
  The matrices $C$ and $\bar C$ are positive semidefinite,
  thus $\lambda_1(C)\ge 0$ and $\lambda_1(\bar C)\ge 0$. By
  $C=\hat q\hat q^\top+(n-2)I_r-\bar C$ the
  statement again follows from \cref{interlacing_rank_one}.
\end{Proof}

We are going to prove the following bounds on the eigenvalues of $C$.
\begin{Theorem}\label{T:condensedeigs}
  Let $b_1^{q_1}\dots b_r^{q_r}$ be an alternating binary sequence
  with $q_1\ge 2$ specifying
  a threshold graph $G$ on $n=\sum_{k=1}^rq_k$ nodes with block degrees
  $p_k$ and $C$ its condensed signless Laplacian.
  \begin{displaymath}
  (0\le )\quad \lambda_1(C)\le p_{r-b_r} \le \lambda_2(C) \le p_{r-b_r-2} \le \dots
  \le \lambda_{\frac{r-b_r-b_1+1}2}(C) \le p_{1+b_1}
\end{displaymath}
and
  \begin{displaymath}
  p_{2-b_1}-1\le \lambda_{\frac{r-b_r-b_1+1}2+1}(C)\le p_{4-b_1}-1\le\lambda_{\frac{r-b_r-b_1+1}2+2}(C)\le
  \dots \le p_{r-(1-b_r)}-1\le \lambda_r(C).
\end{displaymath}
\end{Theorem}
We break the proof into four intermediate steps.
\begin{Observation}\label{O:ce1}
  For $i=1,\dots,\frac{r-b_r-b_1+1}2$ there holds $\lambda_i(C)\le p_{r-b_r-2i+2}$.
\end{Observation}
\begin{Proof}
  Choosing the $r-i$ vectors
  $e_1,e_2,\dots,e_{r-b_r-2i+1},e_{r-b_r-2i+3},e_{r-b_r-2i+5},\dots,e_{r-1+b_r}$ for the right hand
  side $u$ vectors in \cref{T:CF} restricts $x$ to the coordinates
  $e_{r-b_r-2i+2}$, $e_{r-b_r-2i+4}$, $\dots$, $e_{r-b_r}$ (all these indices $k$
  satisfy $b_k=0$), so the bound is obtained by the
  maximum eigenvalue of the submatrix
  \begin{displaymath}
    \left[
      \begin{array}{cccc}
        p_{r-b_r-2i+2} & 0 & \dots & 0 \\
         0   & p_{r-b_r-2i+4} & \ddots &\vdots \\
        \vdots & \ddots & \ddots & 0 \\
        0  & \dots & 0 & p_{r-b_r}
      \end{array}\right].
  \end{displaymath}
  The statement now follows from \cref{eq:pi} and \cref{T:CF}.
\end{Proof}
\begin{Observation}\label{O:ce2}
  For $i=1,\dots,\frac{r-1+b_r+b_1}2$ there holds $p_{r+b_r+1-2i}-1\le \lambda_{r-i+1}(C)$.
\end{Observation}
\begin{Proof}
  Choosing the $r-i$ vectors
  $e_1$, $e_2$, $\dots$, $e_{r+b_r-2i}$, $e_{r+b_r-2i+2}$,
  $e_{r+b_r-2i+4}$, $\dots$, $e_{r-b_r}$ for the left hand
  side $v$ vectors in \cref{T:CF} restricts $x$ to the coordinates
  $e_{r+b_r-2i+1}$, $e_{r+b_r-2i+3}$, $\dots$, $e_{r+br-1}$ (all these indices $k$
  satisfy $b_k=1$), so the bound is obtained by the
  minimum eigenvalue of the submatrix
  \begin{displaymath}
    \bar q\bar q^\top+\left[
      \begin{array}{cccc}
        p_{r+b_r-2i+1}-1 & 0 & \dots & 0 \\
         0   & p_{r+b_r-2i+3}-1 & \ddots &\vdots \\
        \vdots & \ddots & \ddots & 0 \\
        0  & \dots & 0 & p_{r+b_r-1}-1
      \end{array}\right] \text{ with } \bar q=
\left[
  \begin{array}{c}
    \sqrt{q_{r+b_r-2i+1}}\\
    \sqrt{q_{r+b_r-2i+3}}\\
    \vdots\\
    \sqrt{q_{r+b_r-1}}
  \end{array}
\right].
  \end{displaymath}
  Because $\bar q\bar q^\top$ is positive semidefinite, the smallest diagonal
  element is certainly a lower bound, thus
  the statement again follows from \cref{eq:pi} and \cref{T:CF}.
\end{Proof}

\begin{Observation}\label{O:ce3}
  For $i=2,\dots,\frac{r-1+b_r+b_1}2$ there holds $\lambda_{r-i+1}(C) \le p_{r+b_r+3-2i}-1 $.
\end{Observation}
\begin{Proof}
Without loss of generality it suffices to consider the case
$b_r=1$. Indeed, in the case $b_r=0$ we may split off the eigenvector
$e_r$ to eigenvalue $p_r=0$ of $C$ and work with the principal
submatrix $C'$ on indices $i=1,\dots,r-1$. This $C'$ is the condensed
matrix of the threshold graph corresponding to the sequence
$b_1^{q_1}\dots b_{r-1}^{q_{r-1}}$ and its eigenvalues coincide
with the remaining ones of $C$. Thus let $b_r=1$ and
$p_r=n-1$ in the following.

For $i=2$ the statement follows directly from \cref{L:condensedCauchy} by
$\lambda_{r-1}(C)\le n-2$.

For $i>2$ first choose the $i-2$ vectors $u_1=e_r,u_2=e_{r-2},\dots,u_{i-2}=e_{r-2(i-3)}$.
Orthogonality with respect to these vectors already restricts $x$ to the submatrix
  \begin{displaymath}
\left[
      \begin{array}{cccc}
\hat C+\sum_{j=1}^{i-2}q_{r+2-2j}I_{r+4-2i} & 0 & \cdots & 0 \\
   0                                 & p_{r+5-2i} & \ddots & \vdots\\
   \vdots  & \ddots & \ddots & 0 \\
    0      &  \cdots & 0  & p_{r-1}
      \end{array}\right],
  \end{displaymath}
  where $\hat C$ is
  the condensed signless Laplacian of the threshold graph $b_1^{q_1}\dots 1^{q_{r+4-2i}}$. By the same
  argument as above \cref{L:condensedCauchy} yields $\lambda_{r+4-2i-1}(\hat
  C)\le \sum_{j=1}^{r+4-2i}q_j-2$ and
  \begin{displaymath}
    \lambda_{r+4-2i-1}(\hat C+\sum_{j=1}^{i-2}q_{r+2-2j}I_{r+4-2i})\le
  \sum_{j=1}^{r+4-2i}q_j-2+\sum_{j=1}^{i-2}q_{r+2-2j}\overset{\eqref{eq:pidef}}=p_{r+4-2i}-1.
  \end{displaymath}
  Now choose $u_{i-1}$ to hold the Perron vector to the largest eigenvalue of
  $\hat C$ on components $1$ to $r+4-2i$ and zero otherwise, then by
  \cref{eq:pi} and \cref{T:CF} the value $p_{r+4-2i}-1$ is an upper bound on $\lambda_{r-i+1}(C)$.
\end{Proof}
\begin{Observation}\label{O:ce4}
  For $i=2,\dots,\frac{r-b_r-b_1+1}2$ there holds $p_{r-b_r-2i+4} \le \lambda_i(C)$.
\end{Observation}
\begin{Proof}
  Consider the condensed signless Laplacian $\bar C$ of the complement graph
  $(1-b_i)^{q_1}\dots(1-b_r)^{q_r}$ with degrees $\bar p_i$ satisfying
  \cref{eq:compldeg}.
  For $i=2,\dots,\frac{r-1+(1-b_r)+(1-b_1)}2$ \cref{O:ce3} proves
  $\lambda_{r-i+1}(\bar C)\le \bar p_{r+(1-b_r)+3-2i}-1$. Therefore
  \cref{L:condensedCauchy} and \cref{eq:compldeg} establish
  $\lambda_{i}(C)\ge n-2-\lambda_{r-i+1}(\bar C)\ge
   n-2 -(n-1-p_{r-b_r+4-2i}-1)=p_{r-b_r+4-2i}$.
\end{Proof}
Observations \ref{O:ce1}--\ref{O:ce4} prove \cref{T:condensedeigs}. Together
with \cref{L:directvecs} the latter gives rise to the main result,
which is more conveniently stated in terms of the degree sequence
$d_n(G) \geq \cdots \geq d_1(G)$ of the graph. Furthermore, if $b_1=1$ and the number $\bar k$ of ones in the binary sequence exceeds $q_1$, it will be
possible to improve the bound on the larger ``central'' eigenvalue
$\lambda_{\frac{r-b_r-b_1+1}2+1}(C)$ by one.
\begin{Theorem}\label{T:main}
Let $G$ be a threshold graph with $n$ vertices represented by the
binary sequence $b_1^{q_1}\dots b_r^{q_r}$. Denote the degree sequence
of $G$ by $d_n(G) \geq \cdots \geq d_1(G)$ and the number of ones in
the binary sequence of $G$ by $\bar{k}=\sum_{i=1}^rb_iq_i$. Then, the eigenvalues of $G$ satisfy
$$
\lambda_n \ge d_n-1 \ge \lambda_{n-1}\ge \dots \ge \lambda_{n+1-\bar{k}} \ge
  d_{n+1-\bar{k}}-1 \ge d_{n-\bar{k}}\ge \lambda_{n-\bar{k}}\ge
  \dots\ge d_1\ge\lambda_1\ge 0.
$$
Furthermore, if $b_1 = 1$ and $\bar k > q_1$, then $\lambda_{n-\bar{k}+q_1} \geq d_{n-\bar{k}+q_1}$.
\end{Theorem}
\begin{Proof}Note that $b_r=0$ results in appending $q_r$ isolated nodes or,
  equivalently in appending zero rows  and columns which does not
  influence the other eigenvalues and eigenvectors. So for simplifying
  notation we assume without loss of generality
  $b_r=0$. Then according to \eqref{eq:pidef} and \eqref{eq:pi} the
  degrees satisfy
  \begin{displaymath}
  \begin{array}{c}
    d_1=\dots=d_{q_r}=p_r<
    d_{q_r+1}=\dots=d_{q_r+q_{r-2}}=p_{r-2}<\dots\le
    d_{n-\bar{k}}=p_{1+b_1}\le p_{2-b_1}-1,\\
    p_{2-b_1}=d_{n-\bar k+1}=\dots=d_{n-\bar k+q_{2-b_1}}<
    p_{4-b_1}=\dots <p_{r-1}=d_{n-q_{r-1}+1}=\dots=d_n.
\end{array}
\end{displaymath}
With this \cref{T:condensedeigs} and \cref{L:directvecs} yield the general
eigenvalue bounds.

Now assume $b_1 = 1$ and $q_1< \bar k$. This implies $q_2 \geq
1$. Then, by \eqref{eq:pidef} and the relation above, $d_{n-\bar k}=p_2=\bar k-q_1$,
$d_{n-\bar k+1}=\dots=d_{n-\bar k+q_1}=\bar k-1$,  and the $(\overline{k}+1) \times (\overline{k}+1)$ matrix
$$
M = \begin{bmatrix}
d_n & 1 & \cdots & 1 & 1 & 1\\
1 & d_{n-1} & \cdots & 1 & 1 & 1\\
\cdots & & & \cdots & &\cdots\\
1 & 1 & \cdots & \overline{k}-1 & 1& 0\\
1 & 1 & \cdots & 1 & \overline{k}-1 & 0\\
1 & 1 & \cdots &  0 & 0 & \overline{k}-q_1
\end{bmatrix}
$$
is a principal submatrix of $Q(G)$ for a suitable index reordering. The matrix
$$
Q(X)= \begin{bmatrix}
\overline{k} & 1 & \cdots & 1 & 1 & 1\\
1 & \overline{k} & \cdots & 1 & 1 & 1\\
\cdots & & & \cdots & &\cdots\\
1 & 1 & \cdots & \overline{k}-1 & 1& 0\\
1 & 1 & \cdots & 1 & \overline{k}-1 & 0\\
1 & 1 & \cdots &  0 & 0 & \overline{k}-q_1
\end{bmatrix}
$$
is the signless Laplacian matrix of the graph $X$ obtained from $K_{\overline{k}+1}$ by removing a copy of $K_{q_1,1}$. Hence, the matrix $Q(X)$ may be written as
$$
Q(X) = Q(K_{\overline{k}+1}) - Q((\overline{k}-q_1)K_1\cup K_{q_1,1}).
$$
The signless Laplacian $Q(K_{q_1,1})$ has an eigenvalue zero  with
eigenvector $(\mathbf{1}^\top_{q_1},-1)^\top$. Therefore the smallest $q_1$ eigenvalues of the matrix $-
Q((\overline{k}-q_1)K_1\cup K_{q_1,1})$ are its only nonzero eigenvalues. Thus, by the
Weyl inequalities of \cref{lemma:weylIneq},
\begin{align*}
\lambda_{q_1+1}(X) &\geq \lambda_{1}(K_{\overline{k}+1}) +
                     \lambda_{q_1+1}(-Q( (\overline{k}-q_1)K_1\cup K_{q_1,1} ))\\
&= \overline{k}-1 + 0 \\
&= \overline{k}-1.
\end{align*}
%For $n-\overline{k} \leq i \leq n$, $d_i \geq \overline{k}-1$, since every $1$-vertex is connected to all other $1$-vertices, so
Now, we notice that $M = Q(X) + F$, where
$$
F = \mathrm{Diag}(d_n - \overline{k}, d_{n-1}  - \overline{k}, \cdots, d_{n-\overline{k}+q_1+1} - \overline{k},0, \cdots,0).
$$
Therefore, again by \cref{lemma:weylIneq} we obtain
\begin{align*}
\lambda_{q_1+1}(M) &\geq \lambda_{q_1+1}(X) + \lambda_{1}(F)\\
&= \lambda_{q_1+1}(X).
\end{align*}
Since $M$ is a principal submatrix of $Q(G)$, we may apply Cauchy's interlacing
\cref{lemma:cauchy} to obtain
$$
\lambda_{n-\overline{k}+q_1}(G) \geq \lambda_{q_1+1}(M) \geq \lambda_{q_1+1}(X) \geq \overline{k}-1.
$$
Now, the proof follows because when $b_1 = 1$ we have
$d_{n-\bar{k}+q_1} = p_1 = \bar{k} - 1$ by  \eqref{eq:pidef}.
\end{Proof}

\section{Discussion}\label{sec:disc}

We first give an example of the interlacing result and an interpretation
based on the Ferrers diagram of a threshold graph.

\begin{Example}\label{ex1}
Let $G$ be the threshold graph with binary sequence $b_1^{q_1}\dots b_8^{q_{8}}=0^21^2010101^3$. We have $n=12$ and $r=8$. Let $\lambda_i$, $i=1,\ldots,12$, be the eigenvalues of $Q(G)$ and let $\gamma_j$, $j=1,\ldots,8$, be the eigenvalues of $C(G)$ (note that $\gamma_i= \lambda_j$, for some $i$ and $j$).
The signless Laplacian spectrum is (approximately) $$\{ 2.46158, 3.50373, 4.49073, 5.68371, 7, 7, 7.84337, 8.68471, 9.49912, 10, 10, 17.83303\}.$$
As theory predicts,  repeated lines within the Ferrers diagram
lead to signless Laplacian eigenvalues of $G$ determined by \cref{L:directvecs}:
\begin{eqnarray*} \lambda_{10}=\lambda_{11}=10 \ \ \text{and} \ \ \lambda_5 =\lambda_6=7.
\end{eqnarray*}
By the bounds in \cref{T:condensedeigs}, the remaining signless
Laplacian eigenvalues of $G$ are those of $C$ and satisfy:
\begin{displaymath}
  \gamma_1\le 3 \le \gamma_2\le 4 \le \gamma_3\le 5 \le
  \gamma_4\le 7\le \gamma_5\le 8\le \gamma_6\le 9 \le
  \gamma_7\le 10\le \gamma_8.
\end{displaymath}
Note, the bounds hold independently of whether the integral values $p_k-b_k$
are eigenvalues of $G$ or not. Geometrically we can illustrate the
bounds in \cref{T:condensedeigs} by the Ferrers diagram --- its rows
of boxes
display the sorted degree sequence --- in the following way.
\begin{figure}[h]\begin{center}
\begin{tikzpicture}[scale=1.5,auto=left,every node/.style={circle,scale=0.6}]
\draw (-0.3,0.6) rectangle (-0.1,0.8);
\draw (0,0.6) rectangle (0.2,0.8);
\draw (0.3,0.6) rectangle (0.5,0.8);
\draw (0.6,0.6) rectangle (0.8,0.8);
\draw (0.9,0.6) rectangle (1.1,0.8);
\draw (1.2,0.6) rectangle (1.4,0.8);
\draw (1.5,0.6) rectangle (1.7,0.8);
\draw (2,0.6) rectangle (2.2,0.8);
\draw (2.3,0.6) rectangle (2.5,0.8);
\draw (2.6,0.6) rectangle (2.8,0.8);
\draw (2.9,0.6) rectangle (3.1,0.8);
\draw (-0.3,0.3) rectangle (-0.1,0.5);
\draw (0,0.3) rectangle (0.2,0.5);
\draw (0.3,0.3) rectangle (0.5,0.5);
\draw (0.6,0.3) rectangle (0.8,0.5);
\draw (0.9,0.3) rectangle (1.1,0.5);
\draw (1.2,0.3) rectangle (1.4,0.5);
\draw (1.5,0.3) rectangle (1.7,0.5);
\draw (2,0.3) rectangle (2.2,0.5);
\draw (2.3,0.3) rectangle (2.5,0.5);
\draw (2.6,0.3) rectangle (2.8,0.5);
\draw (2.9,0.3) rectangle (3.1,0.5);
\draw (-0.3,0) rectangle (-0.1,0.2);
\draw (0,0) rectangle (0.2,0.2);
\draw (0.3,0) rectangle (0.5,0.2);
\draw (0.6,0) rectangle (0.8,0.2);
\draw (0.9,0) rectangle (1.1,0.2);
\draw (1.2,0) rectangle (1.4,0.2);
\draw (1.5,0) rectangle (1.7,0.2);
\draw (2,0) rectangle (2.2,0.2);
\draw (2.3,0) rectangle (2.5,0.2);
\draw (2.6,0) rectangle (2.8,0.2);
\draw (2.9,0) rectangle (3.1,0.2);
\draw (-0.3,-0.3) rectangle (-0.1,-0.1);
\draw (0,-0.3) rectangle (0.2,-0.1);
\draw (0.3,-0.3) rectangle (0.5,-0.1);
\draw (0.6,-0.3) rectangle (0.8,-0.1);
\draw (0.9,-0.3) rectangle (1.1,-0.1);
\draw (1.2,-0.3) rectangle (1.4,-0.1);
\draw (1.5,-0.3) rectangle (1.7,-0.1);
\draw (2,-0.3) rectangle (2.2,-0.1);
\draw (2.3,-0.3) rectangle (2.5,-0.1);
\draw (2.6,-0.3) rectangle (2.8,-0.1);
\draw (-0.3,-0.6) rectangle (-0.1,-0.4);
\draw (0,-0.6) rectangle (0.2,-0.4);
\draw (0.3,-0.6) rectangle (0.5,-0.4);
\draw (0.6,-0.6) rectangle (0.8,-0.4);
\draw (0.9,-0.6) rectangle (1.1,-0.4);
\draw (1.2,-0.6) rectangle (1.4,-0.4);
\draw (1.5,-0.6) rectangle (1.7,-0.4);
\draw (2,-0.6) rectangle (2.2,-0.4);
\draw (2.3,-0.6) rectangle (2.5,-0.4);
\draw (-0.3,-0.9) rectangle (-0.1,-0.7);
\draw (0,-0.9) rectangle (0.2,-0.7);
\draw (0.3,-0.9) rectangle (0.5,-0.7);
\draw (0.6,-0.9) rectangle (0.8,-0.7);
\draw (0.9,-0.9) rectangle (1.1,-0.7);
\draw (1.2,-0.9) rectangle (1.4,-0.7);
\draw (1.5,-0.9) rectangle (1.7,-0.7);
\draw (2,-0.9) rectangle (2.2,-0.7);
\draw (-0.3,-1.2) rectangle (-0.1,-1);
\draw (0,-1.2) rectangle (0.2,-1);
\draw (0.3,-1.2) rectangle (0.5,-1);
\draw (0.6,-1.2) rectangle (0.8,-1);
\draw (0.9,-1.2) rectangle (1.1,-1);
\draw (1.2,-1.2) rectangle (1.4,-1);
\draw (1.5,-1.2) rectangle (1.7,-1);
\draw (2,-1.2) rectangle (2.2,-1);
\draw (-0.3,-1.7) rectangle (-0.1,-1.5);
\draw (0,-1.7) rectangle (0.2,-1.5);
\draw (0.3,-1.7) rectangle (0.5,-1.5);
\draw (0.6,-1.7) rectangle (0.8,-1.5);
\draw (0.9,-1.7) rectangle (1.1,-1.5);
\draw (1.2,-1.7) rectangle (1.4,-1.5);
\draw (1.5,-1.7) rectangle (1.7,-1.5);
\draw (-0.3,-2.2) rectangle (-0.1,-2);
\draw (0,-2.2) rectangle (0.2,-2);
\draw (0.3,-2.2) rectangle (0.5,-2);
\draw (0.6,-2.2) rectangle (0.8,-2);
\draw (0.9,-2.2) rectangle (1.1,-2);
\draw (1.2,-2.2) rectangle (1.4,-2);
\draw (1.5,-2.2) rectangle (1.7,-2);
\draw (-0.3,-2.5) rectangle (-0.1,-2.3);
\draw (0,-2.5) rectangle (0.2,-2.3);
\draw (0.3,-2.5) rectangle (0.5,-2.3);
\draw (0.6,-2.5) rectangle (0.8,-2.3);
\draw (0.9,-2.5) rectangle (1.1,-2.3);
\draw (-0.3,-2.8) rectangle (-0.1,-2.6);
\draw (0,-2.8) rectangle (0.2,-2.6);
\draw (0.3,-2.8) rectangle (0.5,-2.6);
\draw (0.6,-2.8) rectangle (0.8,-2.6);
\draw (-0.3,-3.1) rectangle (-0.1,-2.9);
\draw (0,-3.1) rectangle (0.2,-2.9);
\draw (0.3,-3.1) rectangle (0.5,-2.9);
\draw [thick] (-0.6,0.86) --(-0.15,0.86);
\node[circle] (14) at (-1.3,0.86) {\Large{$\gamma_8=\lambda_{12}$}};
\draw [thick] (-0.6,0.55) --(-0.15,0.55);
\node[circle] (14) at (-1.3,0.55) {\Large{$10=\lambda_{11}$}};
\draw [thick] (-0.6,0.25) --(-0.15,0.25);
\node[circle] (14) at (-1.3,0.25) {\Large{$10=\lambda_{10}$}};
\draw [thick] (-0.6,-0.05) --(-0.15,-0.05);
\node[circle] (14) at (-1.3,-0.05) {\Large{$\gamma_7=\lambda_9$}};
\draw [thick] (-0.6,-0.35) --(-0.15,-0.35);
\node[circle] (14) at (-1.3,-0.35) {\Large{$\gamma_6=\lambda_8$}};
\draw [thick] (-0.6,-0.65) --(-0.15,-0.65);
\node[circle] (14) at (-1.3,-0.65) {\Large{$\gamma_5=\lambda_7$}};
\draw [thick] (-0.6,-0.95) --(-0.15,-0.95);
\node[circle] (14) at (-1.3,-0.95) {\Large{$7=\lambda_6$}};
\draw [thick] (-0.6,-1.85) --(-0.15,-1.85);
\node[circle] (14) at (-1.3,-1.85) {\Large{$7=\lambda_5$}};

\draw [thick] (-0.6,-2.25) --(-0.15,-2.25);
\node[circle] (14) at (-1.3,-2.25) {\Large {$\gamma_4=\lambda_4$}};

\draw [thick] (-0.6,-2.55) --(-0.15,-2.55);
\node[circle] (14) at (-1.3,-2.55) {\Large {$\gamma_3=\lambda_3$}};

\draw [thick] (-0.6,-2.85) --(-0.15,-2.85);
\node[circle] (14) at (-1.3,-2.85) {\Large {$\gamma_2=\lambda_2$}};

\draw [thick] (-0.6,-3.25) --(-0.15,-3.25);
\node[circle] (14) at (-1.3,-3.25) {\Large{$\gamma_1=\lambda_1$}};
\node[circle] (14) at (4,0.65) {\Large{$b_8$}};
\node[circle] (14) at (4,0.35) {\Large{$b_8$}};
\node[circle] (14) at (4,0.05) {\Large{$b_8$}};
\node[circle] (14) at (4,-0.25) {\Large{$b_6$}};
\node[circle] (14) at (4,-0.55) {\Large{$b_4$}};
\node[circle] (14) at (4,-0.85) {\Large{$b_2$}};
\node[circle] (14) at (4,-1.15) {\Large{$b_2$}};
\node[circle] (14) at (4,-1.6) {\Large{$b_1$}};
\node[circle] (14) at (4,-2.1) {\Large{$b_1$}};
\node[circle] (14) at (4,-2.4) {\Large{$b_3$}};
\node[circle] (14) at (4,-2.7) {\Large{$b_5$}};
\node[circle] (14) at (4,-3.0) {\Large{$b_7$}};
\end{tikzpicture}
\end{center}\caption{Illustration of the interlacing on the Ferrers diagram}\label{fig:Ferrers}
\end{figure}
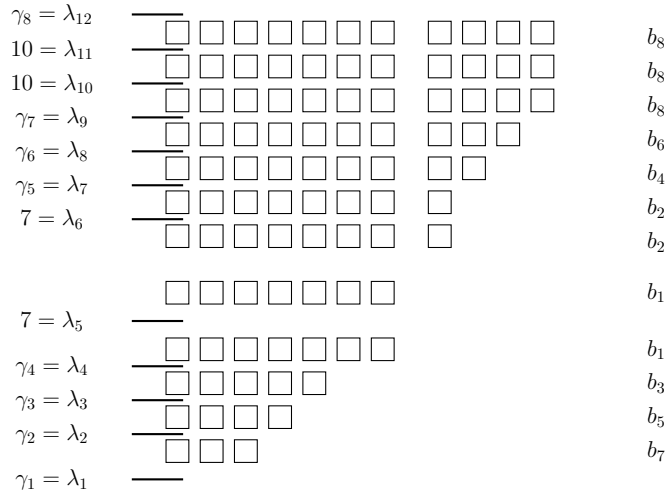
\end{Example}

The next result investigates the development of the eigenvalues upon
appending a one to the binary construction sequence.

\begin{Theorem}\label{T:append_1}
 Let $b_1^{q_1}\dots b_r^{q_r}$ be an alternating binary sequence with
 $q_1\ge 2$ specifying a threshold graph $G$ on $n=\sum_{k=1}^rq_k$ nodes and eigenvalues $0\leq \lambda_1 \leq \lambda_2 \leq \cdots \leq \lambda_n$.
 The eigenvalues $0\leq \lambda^\prime_1 \leq \lambda^\prime_2 \leq \cdots \leq \lambda^\prime_{n+1}$ of the signless Laplacian of the threshold graph $G^\prime$ specified by $b_1^{q_1}\dots b_r^{q_r}1$ satisfy
 \begin{equation}\label{eq:plusone}
   0 \leq \lambda^\prime_i \leq \lambda_i + 1 \leq \lambda^\prime_{i+1} \mbox{~~~~ for } i=1,\ldots,n.
 \end{equation}
Furthermore, $\max\{n + 1,\lambda_n + 2\} \leq\lambda^\prime_{n+1}$.
\end{Theorem}
\begin{Proof} Because
$$ 0 \preceq Q(G^\prime) = \left[
                             \begin{array}{cc}
                               Q(G) & \textbf{0}_{n\times 1}  \\
                               \textbf{0}_{1\times n}& 0 \\
                             \end{array}
                           \right]
+
\sum_{i=1}^n (e_i+e_{n+1})(e_i+e_{n+1})^\top =\left[
                             \begin{array}{cc}
                               Q(G) & \textbf{1}_{n\times 1}  \\
                               \textbf{1}_{1\times n}& n-1 \\
                             \end{array}
                           \right]
+ I_{n+1},$$
inequalities \eqref{eq:plusone} follow directly from Cauchy's
interlacing
\cref{lemma:cauchy}.

$\lambda^\prime_{n+1} \geq n+1$ is a consequence of the maximum eigenvalue of
$\sum_{i=1}^n (e_i+e_{n+1})(e_i+e_{n+1})^\top$
being $n + 1$ (eigenvector $(\textbf{1}_{1\times n},n)^\top$). Suppose
now $\lambda_n+2 > n+1$, \ie, $\lambda_n > n - 1$. Then \cref{T:main}
ensures $ \lambda^\prime_n \leq n-1 < \lambda_n$ and $\trace
Q(G^\prime) = \trace Q(G) + 2n
$ yields $\sum_{i=1}^{n+1} \lambda^\prime_i = \sum_{i=1}^n  \lambda_i + 2n =  \sum_{i=1}^{n-1} (\lambda_i +1)  + \lambda_n  + 1 + n =  \sum_{i=1}^{n-1} (\lambda_i + 1)  + (n-1)+ (\lambda_n  + 2) > \sum_{i=1}^n  \lambda^\prime_i + \lambda_n + 2$, implying that $\lambda^\prime_{n+1} > \lambda_n + 2.$
\end{Proof}

\section{The signless Brouwer conjecture for threshold graphs}\label{sec:brouwer}

For an integer $k$ with $1 \leq k \leq n$, we denote by $S_k(G)$ the sum of
the $k$ largest signless Laplacian eigenvalues of a graph $G$,  that is
$S_k(G)=\sum_{i=n+1-k}^{n} \lambda_i(Q(G))$. Ashraf et al
\cite{ASHRAF20134539} posed the following conjecture.

\begin{Conjecture}[Signless Brouwer Conjecture \cite{ASHRAF20134539}]\label{con1}
Let $G$ be a graph with $n$ vertices. Then
$$S_k(G) \leq |E|+ {k+1 \choose 2},$$
where $k$ is an integer with $1 \leq k \leq n$.
\end{Conjecture}

This conjecture was motivated by Brouwer's Conjecture \cite{bookbrouwer},
which states that the same inequality is valid for the sum of the $k$ largest
Laplacian eigenvalues of a graph $G$ and has been studied by many
researchers. For recent results on Brouwer's Conjecture see
\cite{DU20123672, FRITSCHER2011371, GANIE2016376, HAEMERS20102214,
ROCHA201495, WANG201260,HelmbergTrevisan2017}.

Ashraf el al. \cite{ASHRAF20134539} proposed \cref{con1} and proved it for graphs with at most 10 vertices, for all graphs when $k \in \{1, 2, n-1, n\}$ and for regular graphs. Because trees satisfy Brouwer's Conjecture (see, for example, \cite{FRITSCHER2011371}) and the spectrum of $Q(G)$ is equal to the spectrum of $L(G)$ when $G$ is bipartite, \cref{con1} holds for trees.

In \cite{YANG2014115}, Yang and You studied and proved that \cref{con1} is satisfied by unicyclic graphs, bicyclic graphs and tricyclic graphs with $k \neq 3$. More recently, Chen et al. \cite{Chen2018} proved that \cref{con1} is true for all graphs when $k=n-2$.

In this section, as an application of the main result, we prove that \cref{con1} is true for threshold graphs. We begin by proving some technical lemmas. \cref{degest} provides an expression for $d_{n-\bar{k}}(G)$.

\begin{Lemma}\label{degest}
Let $G$ be a threshold graph with $n$ vertices and alternating binary sequence $b_1^{q_1}\dots b_r^{q_r}$.  We have that $d_{n-\bar{k}}=\bar{k}-q_1b_1$.
\end{Lemma}
\begin{Proof}
It follows from the definition that $\bar{k}=\sum_{h=1}^{r}b_hq_h$ and $d_{n-\bar{k}}$ represents the degree of vertices of the first block of zeros.

If $b_1=0$ then the first block of zeros is the first block of the binary sequence of $G$ and its degree is $p_1=\sum_{h=1}^{n}b_hq_h=\bar{k}-q_1b_1$. If $b_1=1$ then the first block of zeros is the second block of the binary sequence of $G$ and its degree is $p_2=\sum_{h=3}^{r}b_hq_h=\bar{k}-q_1b_1$.
\end{Proof}

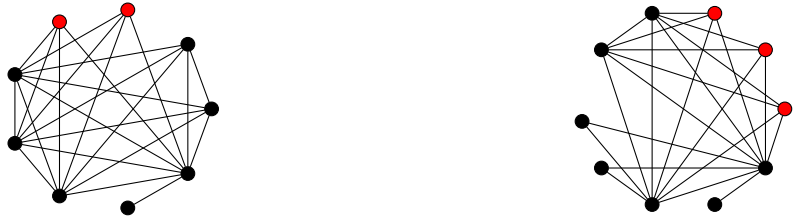
\begin{figure}[H]
\centering
\begin{minipage}{.5\textwidth}
  \centering
\begin{tikzpicture}
\node[draw,circle,fill=black,scale=0.5] (0) at (1.3333,0.0){};
\node[draw,circle,fill=black,scale=0.5] (1) at (1.0213,0.8570){};
\node[draw,circle,fill=red,scale=0.5] (2) at (0.2315,1.3130){};
\node[draw,circle,fill=red,scale=0.5] (3) at (-0.6666,1.1547){};
\node[draw,circle,fill=black,scale=0.5] (4) at (-1.2529,0.4560){};
\node[draw,circle,fill=black,scale=0.5] (5) at (-1.2529,-0.4560){};
\node[draw,circle,fill=black,scale=0.5] (6) at (-0.6666,-1.1547){};
\node[draw,circle,fill=black,scale=0.5] (7) at (0.2315,-1.3130){};
\node[draw,circle,fill=black,scale=0.5] (8) at (1.0213,-0.8570){};
\draw (0) -- (1);
\draw (0) -- (4);
\draw (0) -- (5);
\draw (0) -- (6);
\draw (0) -- (8);
\draw (1) -- (4);
\draw (1) -- (5);
\draw (1) -- (6);
\draw (1) -- (8);
\draw (2) -- (4);
\draw (2) -- (5);
\draw (2) -- (6);
\draw (2) -- (8);
\draw (3) -- (4);
\draw (3) -- (5);
\draw (3) -- (6);
\draw (3) -- (8);
\draw (4) -- (5);
\draw (4) -- (6);
\draw (4) -- (8);
\draw (5) -- (6);
\draw (5) -- (8);
\draw (6) -- (8);
\draw (7) -- (8);
\end{tikzpicture}
\end{minipage}%
\begin{minipage}{.5\textwidth}
  \centering
\begin{tikzpicture}
\node[draw,circle,fill=red,scale=0.5] (0) at (1.3333,0.0){};
\node[draw,circle,fill=red,scale=0.5] (1) at (1.0786,0.7837){};
\node[draw,circle,fill=red,scale=0.5] (2) at (0.4120,1.2680){};
\node[draw,circle,fill=black,scale=0.5] (3) at (-0.4120,1.2680){};
\node[draw,circle,fill=black,scale=0.5] (4) at (-1.0786,0.7837){};
\node[draw,circle,fill=black,scale=0.5] (5) at (-1.3333,-0.1656){};
\node[draw,circle,fill=black,scale=0.5] (6) at (-1.0786,-0.7837){};
\node[draw,circle,fill=black,scale=0.5] (7) at (-0.4120,-1.2680){};
\node[draw,circle,fill=black,scale=0.5] (8) at (0.4120,-1.2680){};
\node[draw,circle,fill=black,scale=0.5] (9) at (1.0786,-0.7837){};
\draw (0) -- (3);
\draw (0) -- (4);
\draw (0) -- (7);
\draw (0) -- (9);
\draw (1) -- (3);
\draw (1) -- (4);
\draw (1) -- (7);
\draw (1) -- (9);
\draw (2) -- (3);
\draw (2) -- (4);
\draw (2) -- (7);
\draw (2) -- (9);
\draw (3) -- (4);
\draw (3) -- (7);
\draw (3) -- (9);
\draw (4) -- (7);
\draw (4) -- (9);
\draw (5) -- (7);
\draw (5) -- (9);
\draw (6) -- (7);
\draw (6) -- (9);
\draw (7) -- (9);
\draw (8) -- (9);
\end{tikzpicture}
\end{minipage}%

\caption{The red vertices have degree $d_{n-\bar k}$.}\label{fig:deg}
\end{figure}
Figure \ref{fig:deg} illustrates \cref{degest}. On the left, we see the threshold graph given by $B_1=1^20^21^301$, in which $n=9, b_1=1, q_1=2$ and $\bar{k}=6$. The red vertices, which represent the first block of zeroes, have degree $p_2=d_{n-\bar{k}}=d_3=\bar{k}-q_1b_1=6-2=4$. On the right, we have the threshold graph given by $B_2=0^31^20^2101$, in which $n=10, b_1=0, q_1=3$ and $\bar{k}=4$. The red vertices (the first block of zeroes) have degree $p_1=d_{n-\bar{k}}=d_6=\bar{k}-q_1b_1=6-2=4$.

In the next lemma we present a lower bound for the number of edges of the graph as a function of its sequence of degrees.

\begin{Lemma}\label{edgeest}
Let $G$ be a threshold graph with $n$ vertices and alternating binary sequence $b_1^{q_1}\dots b_r^{q_r}$. For $k \leq \bar{k}$, we have that

$$\sum_{i=n+1-k}^{n}d_i-{k\choose 2}+{\bar k-k\choose 2}+q_1(\bar k-k)(1-b_1)\le |E|.$$
\end{Lemma}
\begin{Proof}
The term $\sum_{i=n+1-k}^{n}d_i$ adds the degrees of the $k$ last
one-vertices in the binary sequence. In relation to all edges $|E|$
the ${k\choose 2}$ edges between these vertices are counted twice while
at least the ${\bar k-k\choose 2}$ edges among the first $\bar k-k$
one-vertices are not counted in the sum. Furthermore, if $b_1 = 0$, there are $q_1(\bar k-k)$ edges between the last $\bar k - k$ one-vertices and the first $q_1$ zero-vertices that were still not counted. This completes the proof.
\end{Proof}
\begin{figure}[H]
\centering
%\begin{minipage}{.5\textwidth}
%  \centering
\begin{tikzpicture}[scale=2]
\node[draw,circle,fill=black,scale=0.5] (0) at (1.3333,0.0){};
\node[draw,circle,fill=black,scale=0.5] (1) at (1.0213,0.8570){};
\node[draw,circle,fill=black,scale=0.5] (2) at (0.2315,1.3130){};
\node[draw,circle,fill=red,scale=0.5] (3) at (-0.6666,1.1547){};
\node[draw,circle,fill=red,scale=0.5] (4) at (-1.2529,-0.4560){};
\node[draw,circle,fill=black,scale=0.5] (5) at (-1.2529,0.4560){};
\node[draw,circle,fill=black,scale=0.5] (6) at (-0.6666,-1.1547){};
\node[draw,circle,fill=red,scale=0.5] (7) at (0.2315,-1.3130){};
\node[draw,circle,fill=red,scale=0.5] (8) at (1.0213,-0.8570){};
\draw (0) -- (3) [red, very thick];
\draw (0) -- (4) [red, very thick];
\draw (0) -- (7) [black];
\draw (0) -- (8) [black];
\draw (1) -- (3) [red, very thick];
\draw (1) -- (4) [red, very thick];
\draw (1) -- (7) [black];
\draw (1) -- (8)[black];
\draw (2) -- (3) [red, very thick];
\draw (2) -- (4) [red, very thick];
\draw (2) -- (7) [black];
\draw (2) -- (8) [black];
\draw (3) -- (4) [red, very thick];
\draw (3) -- (7) [black];
\draw (3) -- (8)[black];
\draw (4) -- (7) [black];
\draw (4) -- (8) [black];
\draw (4) -- (5)                   [blue, ultra thick];
\draw (5) -- (7) [black] ;
\draw (5) -- (8)[black];
\draw (6) -- (7) [black];
\draw (6) -- (8) [black];
\draw (7) -- (8)[black, ultra thick];
\end{tikzpicture}

\caption{An illustration of the proof of \cref{edgeest}.}\label{fig:lemma19}
\end{figure}
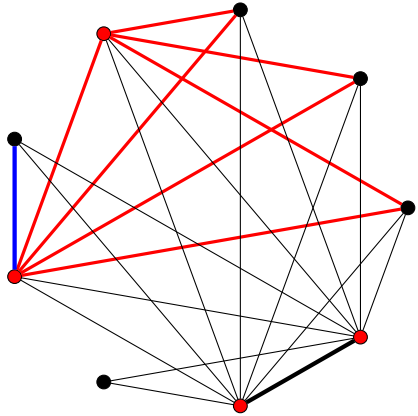
Figure \ref{fig:lemma19} illustrates the proof of \cref{edgeest}. We see the threshold graph given by the binary sequence $B=0^310101^2$. As $n=9$ and $\bar k = 4$, consider $k=2$. The black edges are counted in $d_n + d_{n-1}$, with the thick black edge counted twice. The red edges amount to ${\bar k-k\choose 2}+q_1(\bar k-k)(1-b_1)$. The blue edge is not counted.

Now, we proceed with the proof of the main result of this section.

\begin{Theorem}\label{T:signlessBrouwer}
Let $G$ be a threshold graph with $n$ vertices. For any integer $k$
with $1 \leq k \leq n$ there holds
$$\sum_{i=n+1-k}^{n} \lambda_i(Q(G)) \leq |E|+ {k+1 \choose 2}.$$
\end{Theorem}
\begin{Proof}
Let $b_1^{q_1}\dots b_r^{q_r}$ be the alternating binary sequence with
$q_1\ge 2$ specifying the threshold graph $G$, and $\bar k$ be the number of ones in this sequence. By \cref{T:main},
\begin{equation}\label{eq:lamdeg}
  \lambda_n \ge d_n-1 \ge \dots \ge \lambda_{n+1-\bar{k}}\ge d_{n+1-\bar{k}}-1 \ge d_{n-\bar{k}}\ge \lambda_{n-\bar{k}}\ge \dots\ge d_1\ge\lambda_1.
\end{equation}
Summing the $n+1-\bar{k}$ smallest signless Laplacian eigenvalues we have that
$$\sum_{i=1}^{n-\bar{k}}\lambda_i \geq \sum_{i=1}^{n-1-\bar{k}} d_i.$$
From \cref{degest}, adding $d_{n- \bar{k}}= \bar{k}-q_1 b_1$ on both sides of the above inequality we obtain

$$\sum_{i=1}^{n-\bar{k}}(d_i-\lambda_i) \leq \bar{k}-q_1b_1.$$
The trace $\sum_{i=1}^{n}\lambda_i = \sum_{i=1}^{n} d_i=2|E|$ now
gives rise to
$$\sum_{i=n+1-\bar{k}}^{n}(\lambda_i-d_i) = \sum_{i=1}^{n-\bar{k}}(d_i-\lambda_i) \leq \bar{k}-q_1b_1.$$
Adding $\bar{k}$ on both sides, we obtain
\begin{equation}\label{eq1}
\sum_{i=n+1-\bar{k}}^{n}[\lambda_i-(d_i-1)]  \leq 2\bar{k}-q_1b_1,
\end{equation}
where $\lambda_i-(d_i-1)\ge 0$ for $i=n+1-\bar{k},\dots, n$ by \eqref{eq:lamdeg}.

We first consider $\sum_{i=n+1-k}^{n}\lambda_i$ for $k\le\bar k$. Relation
\eqref{eq1} asserts
\begin{equation}\label{eq2}
\sum_{i=n+1-k}^{n}\lambda_i \leq \sum_{i=n+1-k}^{n}d_i + 2\bar{k}-q_1b_1 - k.
\end{equation}
With \cref{edgeest} this allows one to bound the sum of the largest $k$ eigenvalues by
\begin{align}
  &\sum_{i=n+1-k}^{n}\lambda_i \leq\nonumber\\
  &\qquad\leq |E|+\frac12[4\bar
  k-2q_1b_1-2k+k(k-1)-(\bar k-k)(\bar k-k-1)-2q_1(\bar
  k-k)(1-b_1)] \nonumber\\
  &\qquad= |E|+\frac12[k^2+\bar k -(\bar{k}-k)^2+(\bar{k}-k)(4-2q_1+2q_1b_1) -2q_1b_1]. \label{eig1}
\end{align}
It remains to check that for $k\le \bar k$ the term added to $|E|$ is at
most ${k+1\choose 2}$, or
\begin{displaymath}
  k^2+\bar k -(\bar{k}-k)^2+(\bar{k}-k)(4-2q_1+2q_1b_1) -2q_1b_1\le (k+1)k
\end{displaymath}
which simplifies to
\begin{displaymath}
  -(\bar k-k)^2+(\bar k-k)[5-2q_1(1-b_1)]-2q_1b_1 \leq 0.
\end{displaymath}

Consider first the case $b_1 = 0$. We have $(\bar k -k)[(5-2q_1)-(\bar k -k)] \leq 0$. If $\bar k -k = 0$ the result follows. If $\bar k -k \geq 1$, we notice that $5-2q_1 \leq 1$, since $q_1 \geq 2$, and the result follows.

For $b_1=1$, we have $-(\bar k-k)^2+5(\bar k-k)-2q_1 \leq 0$. We notice that this expression represents a parabola in $(\bar k-k)$, and its maximum ensures that
\begin{displaymath}
-(\bar k-k)^2+5(\bar k-k)-2q_1 \leq \frac{25-8q_1}{4}.
\end{displaymath}
If $q_1 \geq 3$, it follows that $-(\bar k-k)^2+5(\bar k-k)-2q_1 \leq 0$ for $\bar k-k$ integer. If $q_1 =2$, it follows that $-(\bar k-k)^2+5(\bar k-k)-2q_1 \leq 0$ for $\bar k-k \leq 1$ and $\bar k-k \geq 4$.

It remains to consider the cases with $b_1=1$, $q_1 =2$ and $k=\bar k -2, \bar k -3$. We have
$$
\lambda_{n-\overline{k}+2}(G) \geq \overline{k}-1
$$
by \cref{T:main}. Thus, when $k=\overline{k}-2$,
\begin{align*}
\sum_{i=n-\overline{k}+3}^n \lambda_i &= \sum_{i=n-\overline{k}+2}^n \lambda_i - \lambda_{n-\overline{k}+2}(G) \\
&\leq |E| + {\overline{k} \choose 2} - \lambda_{n-\overline{k}+2}(G) \\
&\leq |E| + {\overline{k} \choose 2} - (\overline{k}-1) \\
&= |E| + {\overline{k}-1 \choose 2},
\end{align*}
because the case $k = \bar k - 1$ was already proved. Similarly, when $k = \overline{k}-3$, since $\lambda_{n-\overline{k}+3} \geq \lambda_{n-\overline{k}+2} \geq \overline{k}-1$,
\begin{align*}
\sum_{i=n-\overline{k}+4}^n \lambda_i &= \sum_{i=n-\overline{k}+3}^n \lambda_i - \lambda_{n-\overline{k}+3}(G) \\
&\leq |E|  + {\overline{k}-1 \choose 2} - \lambda_{n-\overline{k}+3}(G) \\
&\leq |E|  + {\overline{k}-1 \choose 2} - (\overline{k}-1) \\
&\leq |E|  + {\overline{k}-2 \choose 2}.
\end{align*}
%
%
%
%%%

For $k \geq \bar{k}+1$ we proceed by induction on $k$. For the
induction basis consider $k = \bar{k}+1$. By $\lambda_{n-\bar{k}} \leq
d_{n-\bar{k}} = \bar{k}-q_1b_1 \leq \bar{k}$ there holds
$$\sum_{i=n-\bar{k}}^{n}\lambda_i = \sum_{i=n+1-\bar{k}}^{n} \lambda_i + \lambda_{n-\bar{k}} \leq |E|+{\bar{k}+1\choose 2}+ \bar{k} \leq |E|+{\bar{k}+2\choose 2}.$$
For  $k > \bar{k}+1$, we have $\lambda_{n+1-k} \leq d_{n+1-k} \leq d_{n-\bar{k}}=\bar{k}-q_1b_1 \leq \bar{k} \leq k$. Consequently, 

$$\sum_{i=n+1-k}^{n}\lambda_i = \sum_{i=n+2-k}^{n} \lambda_i + \lambda_{n+1-k} \leq |E|+{k\choose 2}+ k = |E|+{k+1\choose 2}.$$
This completes the proof.
\end{Proof}

\section*{Acknowledgments}
The authors acknowledge the partial support of DAAD PROBRAL Grant
56267227 - Germany and MATH-AMSUD under project GSA, brazilian team
financed by CAPES under project 88881.694479/2022-01. V. Trevisan also
acknowledges partial support of CNPq grant 310827/2020-5, and FAPERGS grant PqG 17/2551-0001. This research is
part of the doctoral studies of Guilherme Torres and financed in part
by the Coordenação de Aperfeiçoamento de Pessoal de Nível Superior -
Brasil (CAPES) - Finance Code 001.
\bibliographystyle{abbrv}
%\bibliography{biblio190925}

\end{document}